\input amstex

\documentstyle{amsppt}
\NoRunningHeads \magnification=1200 \pagewidth{32pc}
\pageheight{42pc} \vcorrection{1.2pc}

\define\wh{\widehat}
\define\wt{\widetilde}
\define\bc{\Bbb C}
\define\bz{\Bbb Z}

\define\ep{\epsilon}

\topmatter
\title  A quantized Tits-Kantor-Koecher algebra
\endtitle

\footnotetext"$^{1}$"{The author gratefully acknowledges the support
of the Natural Sciences and Engineering Research Council of Canada
and the Chinese Academy of Sciences.} \footnotetext"$^{2}$"{The
author gratefully acknowledges the support from  NSF, NSA and NSFC's
Overseas Distinguished Youth Grant.}

\author Yun Gao
\footnotemark"$^{1}$" and Naihuan Jing \footnotemark"$^{2}$"
\endauthor

\address
YG: Department of Mathematics and Statistics, York University,
Toronto,
\newline Canada M3J 1P3
\endaddress
\email ygao\@yorku.ca
\endemail
\address
NJ: Department of Mathematics, North Carolina State University,
Raleigh, NC 27695-8205, USA
\endaddress
\email jing\@math.ncsu.edu
\endemail

\abstract

 We propose a quantum analogue of a Tits-Kantor-Koecher algebra with a Jordan torus
 as an coordinated
 algebra  by
looking at the vertex operator construction over a Fock space.

\endabstract

\endtopmatter
\document

Quantum toroidal algebras were first introduced by Ginzburg,
Kapranov and Vasserot
  [GKV] in the study of the Langlands reciprocity for algebraic
surfaces. These algebras are quantized analogues for toroidal Lie
algebras of Moody-Rao-Yokonuma [MRY]. Representations of quantum
toroidal algebras have been studied by Varagnolo-Vasserot [VV],
Saito-Takemura-Uglov [STU], Saito [S], Frenkel-Jing-Wang [FJW],
Take\-mura-Uglov[TU], Gao-Jing [GJ1,2], and among others.

The Tits-Kantor-Koecher (TKK) algebra was originally defined from Jordan algebra
in constructing the finite dimensional simple Lie algebras of the exceptional types $E_6$ and $E_7$.
It has also played an important role in
the structure theory of  newly developed extended affine Lie algebras.

\medskip

A  TKK algebra in the extended affine Lie algebras
of type $A_1$ has been realized by gluing a Clifford module and a
Heisenberg module in Tan's paper [T]. This algebra appears as the
core of extended affine Lie algebras of type $A_1$[AABGP] and has
been studied by Yoshii [Y].

\medskip

In this note, we shall propose a quantum analogue of the above
Tits-Kantor-Koecher algebra. Our motivation comes from the vertex
operator construction as was done in [GJ1, GJ2]. We hope that the
quantum TKK algebra will be useful in the study of quantum
toroidal algebras.

Like the quantum Kac-Moody algebra case \cite{J2}
our construction relies on an interesting
combinatorial identity of Hall-Littlewood type \cite{M}.
It suggests that representations of our quantum TKK algebra
probably will provide more generalized combinatorial identities
of this type.
We hope that generalization of our quantum TKK algebras can lead to
further interesting combinatorial structures.

\bigskip

\subhead I. Tits-Kantor-Koecher construction \endsubhead

\medskip

Recall that a Jordan algebra $J$ is a unital commutative algebra
over $\Bbb F$ satisfying
$$ (ab)a^2 = a(ba^2), \text{ for  all } a, b\in J.$$

Note that $J$ may not be associative.

\medskip

\noindent{\bf Example 1.1} Let $\Cal A$ be a unital commutative
associative algebra over $\Bbb F$ and $V$ be an $\Cal A$-module
equipped with an $\Cal A$-bilinear form
$$f: V\times V \to \Cal A.$$
Then
$$J(\Cal{A}, V, f) = \Cal{A}\oplus V$$ becomes a Jordan algebra
over $\Bbb F$ under the product
$$(a + u)(b+ v) =(ab + f(u, v)) +
(av + bu)$$
for $ a, b\in A, u, v\in V.$

\medskip

 Let $J$ be a Jordan algebra. Set
 $$D_{a, b} = [L_a, L_b],$$ for $a,
 b\in J$. The $\Bbb F$-linear span $D_{J, J}$ of all $D_{a, b}$'s is a Lie
 algebra called the inner derivation algebra of $J$. They satisfy
 the following relations:
 $$\align & D_{a, b} + D_{b, a} =0, \\
 & D_{ab, c} + D_{bc, a} + D_{ca, b} = 0,\\
 &[D, D_{a, b}] = D_{Da, b} + D_{a, Db},
 \endalign$$
 for $a, b, c\in J$ and any derivation $D$ of $J$.

\medskip

 The Tits-Kantor-Koecher algebra $K(J)$ is defined to be a Lie algebra
$$K(J)= (sl_2(\Bbb{F})\otimes_\Bbb{F} J) \oplus D_{J, J}$$
with Lie bracket:
$$\align &[A\otimes a, B\otimes b] = [A, B]\otimes ab + 2tr(AB)D_{a, b}, \\\
&[D, A\otimes a]= A\otimes Da,\endalign$$ for $A, B\in
sl_2(\Bbb{F})$, $a, b\in J, D\in D_{J, J}$.

\medskip

In the above example, we let
$$\Bbb{F}=\bc, \Cal{A}=\bc[t_1^{\pm
2}, t_2^{\pm 2}], V=\Cal{A}w_1 \oplus \Cal{A}w_2, f(w_i,
w_j)=\delta_{ij}t_i^2.$$ Let $J$ be the resulting Jordan algebra.
The Tits-Kantor-Koecher algebra $K(J)$ is called a Baby TKK in
[T]. This TKK algebra is indeed the smallest possible core of the extended affine Lie algebra which is coordinated
by a Jordan torus--a nonassociative algebra.

\medskip

Let $d_1, d_2$ be the degree derivations of $J$. Define $\chi: J
\to \bc$ to be the $\bc$-linear function given by
$$\chi(w_1^{n_1}w_2^{n_2}) =\cases 1, & \text{  if  } n_1=n_2=0,\\
0, & \text{ otherwise.}
\endcases$$

\medskip

Define a two dimensional central extension of $K(J)$ as follows:
$$\wh{K}(J)=K(J)\oplus\bc c_1\oplus\bc c_2$$
with Lie bracket
$$\align  [A\otimes a, B\otimes b]
& =[A, B]\otimes ab  + 2tr(AB)D_{a, b}\\
& + tr(AB)\chi((d_1a)b)c_1 + tr(AB)\chi((d_2a)b)c_2 \endalign
$$ where $A, B\in sl_2(\bc), a,b \in J$, and $c_1, c_2$ are
central elements of $\wh{K}(J)$.

\medskip

The semi-direct product of the Lie algebra $\wt{K}(J)$ and the two
degree derivations:
$$\wt{K}(J) =\wh{K}(J)\oplus\bc d_1 \oplus
 \bc d_2$$ is an extended affine Lie algebra which is
 the smallest extended affine Lie algebra beyond the finite and
 affine types.

\medskip

\noindent{\bf Remark 1.2} Note that the Lie algebra $\wh{K}(J)$ is
generated by
$$\align & e_{12}\otimes w_1^m, \quad e_{21}\otimes w_1^m,\quad (e_{11}-e_{22})\otimes w_1^m,
\\
& e_{12}\otimes (w_2w_1^{2m}), \quad  e_{21}\otimes(w_2w_1^{2m}),
\quad (e_{11}-e_{22})\otimes (w_2w_1^{2m})
\endalign$$
for $m\in \Bbb Z$, where $e_{ij}$'s are the standard matrix units.

\bigskip

\subhead  II. A Quantum TKK algebra
\endsubhead

\medskip

The quantum TKK algebra $U_q(\wh{K}(J))$ is the unital associative
algebra  generated by
$$ q^{\pm c/2}, k_1^{\pm}, k_0^{\pm}, x^{\pm}_{1,m}, \psi^{\pm}_{1,
m}, x^{\pm}_{0,2m}, \psi^{\pm}_{0, 2m}, m\in\Bbb Z $$ subject to
the following relations that $q^{\pm c/2}$ is central and
$$\align
&[h_{im}, h_{in}]=\frac{[2m]}{m}[mc]\delta_{m, -n}, \tag 2.1 \\
&[h_{1m}, h_{0n}]=-\frac{[m]}{m}[mc](d^n + d^{-n})\delta_{m, -n}, \tag 2.2\\
&[h_{im}, x^{\pm}_{i,n}]=\pm \frac{[2m]}{m}q^{\mp |m|c/2}x^{\pm}_{i,m+n}, \tag 2.3 \\
&[h_{im}, x^{\pm}_{j,n}]=\mp \frac{[m]}{m}q^{\mp |m|c/2}(d^m + d^{-m})x^{\pm}_{j,m+n}, \tag 2.4\\
&x^{\pm}_{1,m+1}x^{\pm}_{1,n}-q^{\pm
2}x^{\pm}_{1,n}x^{\pm}_{1,m+1}=q^{\pm 2}x^{\pm}_{1,m}x^{\pm}_{1,n+1}-x^{\pm}_{1,n+1}x^{\pm}_{1,m},
\tag 2.5\\
& x^{\pm}_{1,m}x^{\pm}_{0,n+2}+q^{\pm
2}x^{\pm}_{1,m+2}x^{\pm}_{0,n} +x^{\pm}_{0,n}x^{\pm}_{1,m+2}+q^{\pm 2}x^{\pm}_{0,n+2}x^{\pm}_{1,m}=0,
\tag 2.6\\
&[x^+_{im}, x^-_{jn}]=\frac{\delta_{ij}}{q-q^{-1}} (\psi^+_{i,
m+n}q^{(m-n)c/2}-
\psi^-_{i, m+n}q^{(n-m)c/2}), \tag 2.7\\
&x^{\pm}_{i,m_1}x^{\pm}_{i,m_2}x^{\pm}_{i, m_3}x^\pm_{j, n} +
[3]x^{\pm}_{i,m_1}x^{\pm}_{i,m_2}x^{\pm}_{j, n}x^{\pm}_{i,m_3} \\
 & +[3]x^{\pm}_{i,m_1}x^{\pm}_{j, n}x^{\pm}_{i,m_2}x^{\pm}_{i,m_3}+
 x^{\pm}_{j,n}x^{\pm}_{i,m_1}x^{\pm}_{i, m_2}x^\pm_{i, m_3} \tag 2.8\\
 &+ \text{ Perm}\{ m_1, m_2, m_3\}
 =0, \text{ for } i\neq j,
\endalign
$$
where $d=-\sqrt{-1}$,
$$[m]=\frac{q^m-q^{-m}}{q-q^{-1}}, \quad [mc]=\frac{q^{mc}-q^{-mc}}{q-q^{-1}},$$
and
$$\align &\sum_{n=0}^\infty
\psi^{\pm}_{1,n}z^{\mp n}=k_{1}^{\pm }exp(\pm (q-q^{-1})\sum_{n>0}
h_{1, \pm n}z^{\mp n}), \tag 2.9\\
& \sum_{n=0}^\infty \psi^{\pm}_{0,2n}z^{\mp 2n} =k_{0}^{\pm
}exp(\pm (q-q^{-1})\sum_{n>0} h_{0, \pm 2n}z^{\mp 2n}) \tag 2.10
\endalign
$$

Write
$$\align
e_1(z)&=\sum_{n\in\Bbb Z}x^+_{1m}z^{-m}, \quad
f_1(z)=\sum_{n\in\Bbb Z}x^+_{1m}z^{-m}\\
e_0(z)&=\sum_{n\in\Bbb Z}x^+_{02m}z^{-m}, \quad
f_0(z)=\sum_{n\in\Bbb Z}x^+_{02m}z^{-m}
\endalign
$$
Then the defining relations for the quantum TKK algebra can be rewritten
as
$$\align
\psi_i^+(z)\psi_j^-(w)&=\psi_j^-(w)\psi_i^+(z)\frac{q^2c^{-4}(\frac wz)^2+1}
{c^{-4}(\frac wz)^2+q^2}\frac{c^{4}(\frac wz)^2+q^2}{q^2c^{4}(\frac wz)^2+1}, \quad i\neq j
\tag2.11\\
\psi_i^+(z)\psi_i^-(w)&=\psi_i^-(w)\psi_i^+(z)\frac{q^{-2}c^{-2}(\frac wz)-1}
{c^{-2}(\frac wz)-q^{-2}}\frac{c^{2}(\frac wz)-q^{-2}}{q^{-2}c^{2}(\frac wz)-1}\tag2.12\\
(z-q^2w)&e_i(z)e_i(w)=(q^2z-w)e_i(w)e_i(z), \tag 2.13\\
[e_i(z), f_j(w)]&=\frac{\delta_{ij}}{q-q^{-1}}\{\psi_i^+(cw)\delta(c^{-2}\frac zw)
-\psi_i^-(cz)\delta(c^2\frac zw)\} \tag2.14\\
(w^2+q^{-2}z^2)&e_1(z)e_0(w)=-(z^2+q^{-2}w^2)e_0(w)e_1(z) \tag2.15\\
(w^2+q^{2}z^2)&f_1(z)f_0(w)=-(z^2+q^{2}w^2)f_0(w)f_1(z) \tag2.16\\
\psi_i^{\pm}(z)e_i(w)&=e_i(w)\psi_i^{\pm}(z)\frac{q^{\mp 2}c^{-1}(\frac wz)^{\pm 1}-1}
{c^{-1}(\frac wz)^{\pm 1}-q^{\mp 2}} \tag 2.17\\
\psi_i^{\pm}(z)f_i(w)&=f_i(w)\psi_i^{\pm}(z)\frac{q^{\pm 2}c(\frac wz)^{\pm 1}-1}
{c(\frac wz)^{\pm 1}-q^{\pm 2}} \tag 2.18\\
\psi_i^{\pm}(z)e_j(w)&=e_j(w)\psi_i^{\pm}(z)\frac{q^{\mp 2}c^{-2}(\frac wz)^{\pm 2}+1}
{c^{-2}(\frac wz)^{\pm 2}+q^{\mp 2}}, \qquad i\neq j \tag 2.19
\endalign
$$
$$\align
Sym_{z_1, z_2, z_3}&\{e_i(z_1)e_i(z_2)e_i(z_3)e_j(w)+[3]e_i(z_1)e_i(z_2)e_j(w)e_i(z_3)+\\
+[3]e_i&(z_1)e_j(w)e_i(z_2)e_i(z_3)+e_j(w)e_i(z_1)e_i(z_2)e_i(z_3)\}=0, \quad
\text{for}\quad a_{ij}=-2\tag2.20\\
Sym_{z_1, z_2, z_3}&\{f_i(z_1)f_i(z_2)f_i(z_3)f_j(w)+[3]f_i(z_1)f_i(z_2)f_j(w)f_i(z_3)+\\
+[3]f_i&(z_1)f_j(w)f_i(z_2)f_i(z_3)+f_j(w)f_i(z_1)f_i(z_2)f_i(z_3)\}=0, \quad
\text{for}\quad a_{ij}=-2\tag2.21
\endalign
$$

 \noindent{\bf Remarks 2.22}

1. The subalgebra generated by
$h_{1m}, x^{\pm}_{1m}$ or $\phi_{im}^{\pm}, x_{1m}^{\pm}$
is isomorphic to the
quantum affine algebra $U_q(\hat{sl}_2)$.

2. The deformation $U_q(\wh{K}(J))$ to
$U(\wh{K}(J))$ can be achieved via

$$\align & x^+_{1, m} \to e_{12}\otimes w_1^m, \quad x^-_{1, m} \to e_{21}\otimes w_1^m, \\
& h_{1m}\to (e_{11}-e_{22})\otimes w_1^m,
\\
& x^+_{0, 2m} \to e_{12}\otimes(w_2w_1^{2m}), \quad x^-_{0, 2m}
\to
e_{21}\otimes(w_2w_1^{2m}),\\
& h_{02m} \to (e_{11}-e_{22})\otimes (w_2w_1^{2m}),\\
&q^{c/2} \to c_1, \quad  k_0^+k_1^+ \to c_2
\endalign$$
for $m\in \Bbb Z$.

\bigskip

\subhead III. Vertex operator representation \endsubhead

\medskip

Let $P=\bz\ep_1\oplus\bz\ep_{2}$ be a rank $2$ free abelian group
provided with a $\bz$-bilinear form $(\cdot ,\cdot)$ defined by
$(\ep_i, \ep_j) = \delta_{ij}$, $1\leq i, j\leq 2$. Let
$Q=\bz(\ep_1-\ep_2)$ be the rank $1$ free subgroup of $P$.

\bigskip

Let
$$\bc[Q] =\sum \oplus\bc e^\alpha $$ be the group algebra of
$Q$. Also, for $\beta\in H = Q\otimes_{\Bbb Z}\Bbb C$, define
$\beta(0)\in \text{End}\bc[Q]$ by
$$\beta(0) e^\alpha = (\beta, \alpha) e^\alpha, \, \text{ for } \alpha\in Q.$$

\medskip

Next let $\ep_i(n)$ and $C$ be the generators of the Heisenberg
algebra $\Cal H$, $ 1\leq i \leq 2, n\in\bz\setminus\{0\}$,
subject to relations that $C$ is central and
$$
[\ep_i(m), \ep_j(n)]=m\delta_{ij}\delta_{m+n, 0}C.\tag 3.1
$$

 \medskip

Let
$$S(\Cal{H}^-) = \bc[\ep_i(n): 1\leq i\leq 2, n\in -\bz_+]$$
denote the symmetric algebra of $\Cal{H}^-$, which is the algebra
of polynomials in infinitely many variables $\ep_i(n), 1\leq i\leq
2, n\in -\bz_+$, where $\bz_+ =\{ n\in \bz: n>0\}$. $S(\Cal{H}^-)$
is
 an $\Cal H$-module in which $C=1$, $\ep_i(n)$ acts as the
multiplication operator for $n\in -\bz_+$, and $\ep_i(n)$ acts as the partial
 differential operator for $n\in \bz_+$.

\medskip

Set
$$V_Q = S(\Cal{H}^-)\otimes \bc[Q].$$
The operator $z^\alpha \in (\text{End}\bc[Q])[z, z^{-1}]$ is
defined as
$$z^\alpha e^\beta = z^{(\alpha, \beta)}e^\beta$$ for
$\alpha, \beta\in Q$.

\medskip

Let $\mu$ be any non-zero complex number. Consider the valuation
$\mu^\alpha$ of the operator $z^\alpha$. Namely, $\mu^\alpha$ is
the operator $\bc[Q]\to \bc[Q]$ given by
$$\mu^\alpha e^\beta = \mu^{(\alpha, \beta)} e^\beta, \, \text{ for } \alpha, \beta\in Q.$$

Now we set
$$\ep_{i+2} = \ep_i, \quad \text{ for } i\in \Bbb Z.$$
Accordingly,
$$(\ep_i, \ep_j) = \delta_{ij} = \delta_{\bar{i}, \bar{j}}, \text{ for } \bar{i}, \bar{j}\in \Bbb{Z}/2\Bbb{Z}.$$

 For $r,  i, j\in \Bbb Z$, we define the vertex operator
$X_{ij}(r, z)$ as follows.
$$
\align
X_{ij}(r, z)
=&
:\exp(-\sum_{n\neq 0}\frac{(\ep_{i}(n)-(-1)^{-rn}q^{(i-j)|n|}\ep_{j}(n))}{n}z^{-n}):  \\
&\qquad e^{\ep_i-\ep_j}z^{\ep_i-\ep_j+\frac{(\ep_i-\ep_j, \ep_i
-\ep_j)}{2}}(-1)^{-r\ep_j- \frac{(\ep_j, \ep_i -\ep_j)}{2}r}
\endalign
$$

Due to $\ep_0=\ep_2$ we note that $X_{01}(r, z)=X_{21}(r, z)$.
Next, for $r, i, j\in \Bbb Z$, and $i \neq j $,  we define
$$\align & u_{ij}(r, z)\\
 = &-q^{(j-i)(\ep_i -\ep_j)}\cdot \exp(\sum_{n\geq 1}\frac{q^{(j-i)n}-q^{(i-j)n}}{n}
(q^{\frac{j-i}{2}n}\ep_i(n)-(-1)^{-nr}q^{\frac{i-j}{2}n}\ep_j(n)z^{-n})\\
&\\
 & v_{ij}(r, z) \\
= & -q^{(i-j)(\ep_i -\ep_j)} \cdot \exp( \sum_{n\geq
1}\frac{q^{(i-j)n}-q^{(j-i)n}}{n}(q^{\frac{j-i}{2}n}\ep_i(-n)-(-1)^{nr}q^{\frac{i-j}{2}n}
\ep_j(-n))z^{n}).
\endalign $$

Write
$$\align & X_{ij}(r, z) = \sum_{n\in\Bbb{Z}} X_{ij}(r, n)z^{-n},\\
& u_{ij}(r, z)=\sum_{n=0}^\infty u_{ij}(r, n)z^{-n}, \\
& v_{ij}(r, z)= \sum_{n=0}^\infty v_{ij}(r, n)z^{n}.
\endalign$$

We now state our main result of this note.

\proclaim{Theorem 3.2} The linear map $\pi$ given by
$$\align & \pi(x^+_{1, m}) = X_{12}(0,m), \quad \pi(x^-_{1, m})= X_{21}(0,
m),\\
& \pi(\psi^+_{1, m}) = u_{12}(0, m),\quad \pi(\psi^-_{1, m}) = v_{12}(0, m)\\
& \pi(x^+_{0, 2m})= X_{01}(1, 2m), \quad \pi(x^-_{0,2m}) = X_{10}(-1, 2m),\\
& \pi(\psi^+_{0, 2m}) = u_{01}(1, 2m), \quad \pi(\psi^-_{0, 2m}) = v_{01}(1, 2m) , \\
& \pi(k_1^\pm) = q^{\pm(\ep_1-\ep_2)}, \quad \pi(k_0^\pm) =
q^{\pm(\ep_2-\ep_1)}, \quad \pi(q^{c/2}) = q^{1/2}
\endalign $$
gives  a representation of $U_q(\wh{K}(J))$.
\endproclaim

\noindent{Proof.} To prove the theorem we need to verify that the defined operators satisfy the
commutation relations in the quantum TKK algebra.

The following result is from \cite{GJ2}.

\proclaim{Lemma 3.3} For $r_1, r_2\in\Bbb Z$ we have
$$\align
X_{ij}(r_1, z)X_{ij}(r_2, z)&=:X_{ij}(r_1, z)X_{ij}(r_2, z):
\frac zw(1-\frac wz)(1-(-1)^{r_2-r_1}q^{-2}\frac wz)(-1)^{r_1}\\
X_{ij}(r_1, z)X_{ji}(r_1, z)&=:X_{ij}(r_1, z)X_{ji}(r_1, z):
\frac wz(1-(-1)^{r_2}\frac w{qz})^{-1}(1-(-1)^{r_1}\frac w{qz})^{-1}(-1)^{r_1}
\endalign
$$
\endproclaim

First of all we notice that the operators $X_{12}(0, z), X_{21}(0, z)$ and
$u_{12}(0, z), v_{12}(0, z)$ gives a level one
representation of the quantum affine algebra $U_q(\hat{sl}_2)$
on the space $V$. In fact let $c=q$ we have
$$
[\ep_1(m)-(-1)^{-rm}q^{-|m|}\ep_2(m),\ep_1(m)-(-1)^{-rn}q^{-|n|}\ep_2(n)]
=m(1+q^{-2|n|})\delta_{m, -n}
$$
Using the identity $e^Ae^B=e^Be^Ae^{[A, B]}$ when $[A, B]$ commutes with
$A$ and $B$, one obtains immediately that
$$\align
(z-q^2w)X_{12}(0, z)X_{12}(0, w)&=(q^2-w)X_{12}(0, w)X_{12}(0, z)\\
[X_{12}(0, z),X_{12}(0, w)]&=\frac1{q-q^{-1}}(u_{12}(cw)\delta(c^{-2}\frac zw)
-v_{12}(cz)\delta(c^2\frac zx))
\endalign
$$

Taking derivative on the operator $u_{ij}(r, z)$ and $v_{ij}(r, z)$,
the map $\pi$ in terms of components is given by
$$
\align
h_{0m}&\to (q^{|m|/2}\ep_1(m)-q^{-|m|/2}\ep_2(m))\frac{[m]}m\\
h_{1m}&\to (q^{|m|/2}\ep_2(m)-(-1)^mq^{-|m|/2}\ep_1(m))d^{-m}\frac{[m]}m
\endalign
$$
It follows that
$$\align
[\pi(h_{1m}), \pi(h_{0n})]&=-\frac{[m]^2}{m}(1+d^{2n})d^{-n}\delta_{m, -n}\\
&=-\frac{[m]}{m}(d^m+d^{-m})\delta_{m, -n},
\endalign
$$
with $c=q$ and $C=1$. Similarly one can check that for $i\neq j$ we have
$$
[\pi(h_{im}), \pi(x_{jn}^{\pm})]=\mp\frac{[m]}mq^{\pm|m|/2}(d^{m}+d^{-m})\pi(x_{j, m+n}^{\pm}).
$$

We now prove that $[x_{im}^+, x_{jn}^-]=0$ for $i\neq j$. In fact we have
$$\align
X_i^+(z)X_j^-(w)&=:X_i^+(z)X_j^-(w):\frac zw(1-\frac wz)(1-p^{-1}\frac wz)p\\
X_j^+(w)X_i^-(z)&=:X_j^+(w)X_i^-(z):\frac wz(1-\frac zw)(1-p\frac zw)p
\endalign
$$
It follows quickly that $[X_i^+(z), X_j^-(w)]=0$.

Finally let's prove the Serre relation.

The following OPE's are direct consequences of
Lemma 3.3.
$$\align
E_1(z)E_1(w)&=:E_1(z)E_1(w):\frac{(z-w)(z-q^{-2}w)}{zw}\\
E_1(z)E_0(w)&=:E_1(z)E_0(w):\frac{zwd}{z^2+q^{-2}w^2}\\
E_1(w)E_0(z)&=:E_1(w)E_0(z):\frac{p^{-1}zwd}{w^2+q^{-2}z^2}
\endalign
$$
Then we have
$$\align
E_1(z_1)&E_1(z_2)E_1(z_3)E_0(w)=:E_1(z_1)E_1(z_2)E_1(z_3)E_0(w):\\
&\cdot\prod_{i<j}
\frac{(z_i-z_j)(z_i-q^{-2}z_j)}{z_iz_j}\cdot\prod_{i=1}^3\frac{z_iwd}{z_i^2+q^{-2}w^2}\\
E_1(z_1)&E_1(z_2)E_0(w)E_1(z_3)=:E_1(z)E_1(z_2)E_0(w)E_1(z_3):\\
&\cdot\prod_{i<j}
\frac{(z_i-z_j)(z_i-q^{-2}z_j)}{z_iz_j}\cdot\frac{p^{-1}z_1z_2z_3w^3d^3}
{(z_1^2+q^{-2}w^2)(z_1^2+q^{-2}w^2)(w^2+q^{-2}z_3^2)}\\
E_1(z_1)&E_0(w)E_1(z_2)E_1(z_3)=:E_1(z)E_0(w)E_1(z_2)E_1(z_3):\\
&\cdot\prod_{i<j}
\frac{(z_i-z_j)(z_i-q^{-2}z_j)}{z_iz_j}\cdot\frac{p^{-2}z_1z_2z_3w^3d^3}
{(z_1^2+q^{-2}w^2)(w^2+q^{-2}z_2^2)(w^2+q^{-2}z_3^2)}\\
E_0(w)&E_1(z_1)E_1(z_2)E_1(z_3)=:E_0(w)E_1(z_1)E_1(z_2)E_1(z_3):\\
&\cdot\prod_{i<j}
\frac{(z_i-z_j)(z_i-q^{-2}z_j)}{z_iz_j}\cdot\frac{p^{-3}z_1z_2z_3w^3d^3}
{(w^2+q^{-2}z_1^2)(w^2+q^{-2}z_2^2)(w^2+q^{-2}z_3^2)}
\endalign
$$
Therefore we have
$$\align
&E_1(z_1)E_1(z_2)E_1(z_3)E_0(w)+[3]E_1(z_1)E_1(z_2)E_0(w)E_1(z_3)+\\
&+[3]E_1(z_1)E_0(w)E_1(z_2)E_1(z_3)+E_0(w)E_1(z_1)E_1(z_2)E_1(z_3)\\
&=:E_1(z_1)E_1(z_2)E_1(z_3)E_0(w):\prod_{i<j}
\frac{(z_i-z_j)(z_i-q^{-2}z_j)}{z_iz_j}(z_1z_2z_3)(wd)^3\\
\cdot\{&\frac1{(z_1^2+q^{-2}w^2)(z_1^2+q^{-2}w^2)(z_1^2+q^{-2}w^2)}
-\frac{[3]}{(z_1^2+q^{-2}w^2)(z_1^2+q^{-2}w^2)(w^2+q^{-2}z_3^2)}\\
+&\frac{[3]}{(z_1^2+q^{-2}w^2)(w^2+q^{-2}z_2^2)(w^2+q^{-2}z_3^2)}
-\frac{1}{(w^2+q^{-2}z_1^2)(w^2+q^{-2}z_2^2)(w^2+q^{-2}z_3^2)}\}\\
&=:E_1(z_1)E_1(z_2)E_1(z_3)E_0(w):\frac{(wd)^3}{z_1z_2z_3}
\prod_{i<j}(z_i-z_j)^2\prod_{i=1}^3\frac1{(z_i^2+q^{-2}w^2)(w^2+q^{-2}z_i^2)}\\
&\cdot\{(w^2+q^{-2}z_1^2)(w^2+q^{-2}z_2^2)(w^2+q^{-2}z_3^2)-[3]
(w^2+q^{-2}z_1^2)(w^2+q^{-2}z_2^2)(z_3^2+q^{-2}w^2)\\
&+[3](w^2+q^{-2}z_1^2)(z_2^2+q^{-2}w^2)(z_3^2+q^{-2}w^2)-(z_1^2+q^{-2}w^2)(z_2^2+q^{-2}w^2)(z_3^2+q^{-2}w^2)\}\\
&\cdot\prod_{i<j}\frac{z_i-q^{-2}z_j}{z_i-z_j}
\endalign
$$
Thus the Serre relation holds if the following combinatorial identity
is truth.

\proclaim{Lemma 3.4} Let $\goth S_3$ act on $z_1, z_2, z_3$ via $\sigma.z_i=z_{\sigma(i)}$. Then
$$\align
\sum_{\sigma\in\goth S_3}&\sigma.[(w^2+q^{-2}z_1^2)(w^2+q^{-2}z_2^2)(w^2+q^{-2}z_3^2)-\\
&-[3](w^2+q^{-2}z_1^2)(w^2+q^{-2}z_2^2)(z_3^2+q^{-2}w^2)\\
&+[3](w^2+q^{-2}z_1^2)(z_2^2+q^{-2}w^2)(z_3^2+q^{-2}w^2)\\&-(z_1^2+q^{-2}w^2)(z_2^2+q^{-2}w^2)(z_3^2+q^{-2}w^2)]\prod_{i<j}\frac{z_i-q^{-2}z_j}{z_i-z_j}=0. \tag 3.5
\endalign
$$
\endproclaim

\noindent{Proof of the Lemma}. Considering the left-hand side as a polynomial in $w$, we
extract the constant term.
$$
\sum_{\sigma\in\goth S_3}(q^{-6}-[3]q^{-4}+[3]q^{-2}-1)(z_1z_2z_3)^2
\sigma.\prod_{i<j}\frac{z_i-q^{-2}z_j}{z_i-z_j}=0.
$$
Similarly the highest coefficient of $w^6$ is seen to be zero.

The coefficient of $w^2$ and $w^4$ are essentially the same up to swapping
of $z_i$ with $z_i^{-1}$. Thus the identity (3.5) in Lemma 3.4 boils down to the truth
of the following identity.
$$
\sum_{\sigma\in\goth S_3}\sigma.\{q^{-3}z_1^2-(q+q^{-1})z_2^2
+q^3z_3^2\}\prod_{i<j}\frac{z_i-q^{-2}z_j}{z_i-z_j}=0, \tag 3.6
$$
where the left-hand side times $q^{-5}-q^{-1}$ is the coefficient of
$w^4$ of the polynomial in Eq. (3.5).

The identity (3.6) is easily proved by comparing coefficients of $z_i$ or
direct verification. Hence Lemma 3.4 is proved.
Similarly one can prove the Serre relations for the $F_i(z)$'s,
and Theorem 3.2 is proved.

\enddocument